%% file: symplectickite2.tex
\documentclass[11pt]{article}
\usepackage{graphics}
\usepackage{graphicx}
\usepackage{color}
\usepackage[all]{xy}
\newfont{\bbb}{msbm10 scaled\magstephalf}
\newfont{\sbbb}{msbm7 scaled\magstephalf}

\def\C{\mbox{\bbb{C}}}

\def\R{\mbox{\bbb{R}}}
\def\Z{\mbox{\bbb{Z}}}
\def\cd{\C^d}
\def\rd{\R^d}

\def\rn{\R^n}
\def\rtwo{\R^2}
\def\zd{\Z^d}
\def\td{T^d}
\def\rddu{(\rd)^*}
\def\rtwodu{(\R^2)^*}
\def\rndu{(\rn)^*}

\def\et1{e^{2\pi i\theta_1}}
\def\etd{e^{2\pi i\theta_d}}

\def\vz{\underline{z}}
\def\zjs{|z_j|^2}

\def\cl{{\cal C}_\lambda}
\def\c0{{\cal C}_0}
\def\D{\Delta}
\def\Dsh{\D^{\sharp}}
\def\xd{X_1,\ldots,X_d}
\def\ld{\lambda_1,\ldots,\lambda_d}

\def\G{\Gamma}
\def\Gsh{\G^{\sharp}}

\def\gsh{\gamma^{\sharp}}
\def\vt{\tilde{V}}
\def\ut{\tilde{U}}
\def\fia{\tau_{\alpha}}
\def\ga{\Gamma_{\alpha}}

\def\fib{\tau_{\beta}}
\def\gb{\G_{\beta}}
\def\vb{V_{\beta}}
\def\vtb{\vt_{\beta}}
\def\gab{g_{\alpha\beta}}
\def\ft{\tilde{f}}
\def\wt{\tilde{W}}
\def\vsh{V^{\sharp}}
\def\wsh{W^{\sharp}}
\def\fbar{\bar{f}}
\def\fsh{f^{\sharp}}

\newfont{\frak}{eufm10 scaled\magstep1}
\newfont{\sfrak}{eufm8 scaled\magstep1}

\def\n{\mbox{\frak n}}

\def\uapr{U_{\alpha'}}
\def\fiapr{\tau_{\alpha'}}
\def\utapr{\tilde{U}_{\alpha'}}
\def\gapr{\G_{\alpha'}}

\def\db{\Delta_{\beta}}
\def\dsh{\Delta_{\sharp}}
\def\uz{U_{\zeta}}

\def\ut{\tilde{U}}
\def\uta{\tilde{U}_{\alpha}}
\def\utb{\tilde{U}_{\beta}}

\def\ua{U_{\alpha}}
\def\ub{U_{\beta}}

\def\piabx{\pi_{\alpha\beta,x}}
\def\pibax{\pi_{\beta\alpha,x}}

\def\pa{p_{\alpha}}
\def\pb{p_{\beta}}

\def\gzb{g_{\zeta\beta}}

\def\ush{U^{\sharp}}

\def\kb{k_{\beta}}

\def\gaz{g_{\alpha\zeta}}
\def\gazw{g_{\alpha\zeta,W}}

\def\gzbw{g_{\zeta\beta,W}}
\def\gabw{g_{\alpha\beta,W}}

\def\fiz{\tau_{\zeta}}
\def\gabx{g_{\alpha\beta,x}}

\def\fd{f^{\dagger}}
\def\ud{U^{\dagger}}
\def\vd{V^{\dagger}}
\def\Gd{\G^{\dagger}}
\def\Dd{\D^{\dagger}}

\def\whs{W^{\sharp}}
\newtheorem{thm}{Theorem}[section]
\newtheorem{defn}[thm]{Definition}
\newtheorem{prop}[thm]{Proposition}
\newtheorem{remark}[thm]{Remark}
\newtheorem{lemma}[thm]{Lemma}

\def\squareforqed{\hbox{\rlap{$\sqcap$}$\sqcup$}}
\def\qed{\ifmmode\else\unskip\quad\fi\squareforqed}
\def\smartqed{\def\qed{\ifmmode\squareforqed\else{\unskip\nobreak\hfil
\penalty50\hskip1em\null\nobreak\hfil\squareforqed
\parfillskip=0pt\finalhyphendemerits=0\endgraf}\fi}}

\newcommand{\proof}{\mbox{\bf Proof.\ \ }}

\newcounter{sect}

\title{\sc The Symplectic Penrose Kite}
\author{\sc Fiammetta Battaglia and Elisa Prato}
\date{}
\begin{document}
\maketitle
\let\thefootnote\relax\footnotetext{Research partially supported by MIUR (Geometria Differenziale e Analisi Globale, PRIN 2007).}
\begin{abstract}
The purpose of this article is to view the Penrose kite from the
perspective of symplectic geometry.
\end{abstract}

\section*{Introduction}
The kite in a Penrose nonperiodic tiling by kite and darts
\cite{pen1,pen2} is an example of a simple convex polytope. By the
Atiyah, Guillemin--Sternberg convexity theorem \cite{a, gs}, the
image of the moment mapping for a Hamiltonian torus action on a
compact connected symplectic manifold is a rational convex polytope.
Moreover, the Delzant theorem \cite{d} provides an exact
correspondence between symplectic toric manifolds and  simple convex
rational polytopes that satisfy a special integrality condition; a
crucial feature of this theorem is that it gives an explicit
construction of the manifold that is associated to each polytope.
The Penrose kite, however, is the most elementary and beautiful
example of a simple convex polytope that is {\em not} rational. The
purpose of this article is to apply to the kite a generalization of
the Delzant construction for nonrational polytopes, which was
introduced by the second--named author in \cite{p1}. We recall that
this generalized construction allows to associate to {\em any}
simple convex polytope $\D$ in $(\R^k)^*$ a $2k$--dimensional
compact connected symplectic {\em quasifold}. Quasifolds are a
natural generalization of manifolds and orbifolds introduced in
\cite{p1}: a local $n$--dimensional model is given by the quotient
of an $n$--dimensional manifold by the smooth action of a discrete
group. In the generalized construction the lattice of the rational
case is replaced by a {\em quasilattice} $Q$, which is the
$\Z$--span of a set of generators of $\R^k$. The torus is replaced
accordingly by a {\em quasitorus}, which is the quotient of $\R^k$
modulo $Q$. The action of the quasitorus on the quasifold is smooth,
effective and Hamiltonian and, exactly as in the Delzant case, the
image of the corresponding moment mapping is the polytope $\D$.

In order to apply the generalized Delzant construction to the kite
we need to choose a suitable quasilattice $Q$, and a set of four
vectors in $Q$ that are orthogonal to the edges of the kite and that
point towards the interior of the polytope. The most natural choice
is to consider, among the various inward--pointing orthogonal
vectors, those four which have the same length as the longest edge
of the kite, and then to choose $Q$ to be the quasilattice that they
generate. We remark that these choices are justified by the geometry
of the kite, and, more globally, by the geometry of any kite and
dart tiling, in the following sense. Let us consider the
quasilattice $R$ which is generated by the vectors that are dual to
the generators of $Q$; notice that the generators of $R$ are {\em
parallel} to the edges of the kite. Then the quasilattice $R$
contains the four vertices of the kite. Moreover, given any kite and
dart tiling, if we consider one of its kites and the associated
quasilattice $R$, then {\em all} of the vertices of the tiling lie
in $R$ and the Delzant procedure can be applied, with respect to
$R$, to each kite of the tiling, giving rise to a unique symplectic
quasifold $M^+_0$ (Theorem~\ref{uguali}).

The four--dimensional quasifold $M^+_0$ turns out to be a very nice
example of a quasifold that is {\em not} a global quotient of a
manifold modulo the smooth action of a discrete group
(Theorem~\ref{noneglobale}), as is the case instead with the
symplectic quasifolds that have been associated to a Penrose rhombus
tiling in \cite{rhombus}.

Quasilattices arise naturally also in the study of the physics of
{\em quasicrystals}. Quasicrystals are some very special alloys,
which were discovered  by Shechtman, Blech, Gratias and Cahn in 1982
\cite{quasicristalli}, that have discrete but non--periodic
diffraction patterns. We remark that the quasilattice $R$ describes
the diffraction pattern of quasicrystals with pentagonal axial
symmetry, which is prohibited for ordinary crystals
\cite{quasicristallireview}. Another symmetry that is forbidden for
crystals but is allowed for quasicrystals is icosahedral symmetry.
In this case too there is a quasilattice underlying both the
structure of the quasicrystals and the suitable analogues of Penrose
tilings in dimension $3$. A symplectic interpretation of this case
is given in \cite{3d}.

The paper is structured as follows. In Section~1 we recall the
classical construction of the Penrose kite and dart from the
pentagram. In Section~2 we introduce the quasilattices $Q$ and $R$
and we discuss their relevant properties. In Section~3 we sketch the
generalized Delzant procedure. In Section~4 we apply the procedure
to the kite and we also describe part of an atlas for the
corresponding quasifold. In Section~5 we show that each kite of any
kite and dart tiling yields the same symplectic quasifold. In
Section~6 we prove that this quasifold is not a global quotient of a
manifold modulo the smooth action of a discrete group. Finally in
the Appendix we recall and partly reformulate the definitions of
quasifold and of related geometrical objects.

\section{The Penrose Kite and Dart}\label{figurekiteedart}

\begin{figure}
\begin{center}
\includegraphics{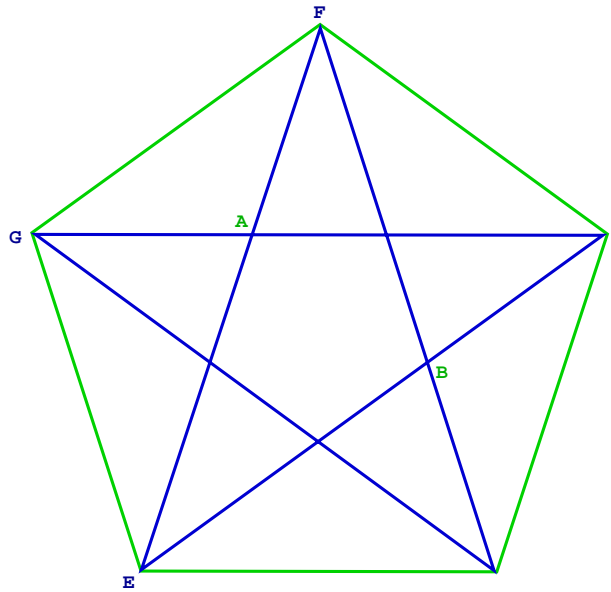}
\end{center}
\caption{The pentagram} \label{pentagramforkite}
\end{figure}
\begin{figure}
\begin{center}
\includegraphics{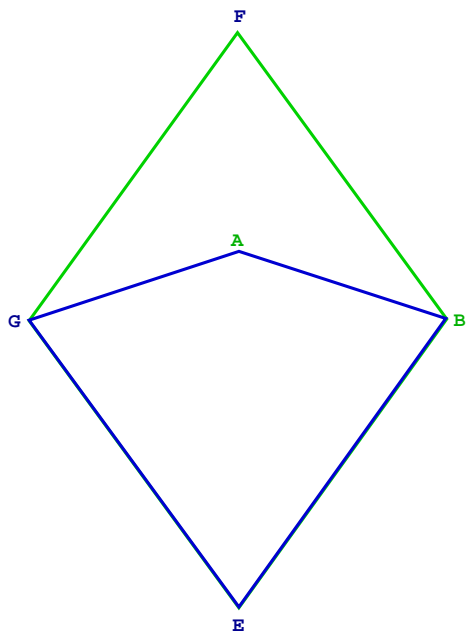}
\end{center}
\caption{The Penrose kite, dart and thick rhombus} \label{kite}
\end{figure}

Let us now recall the procedure for obtaining the Penrose kite and
dart from the pentagram. For a proof of the facts that are needed we
refer the reader to \cite{l}, and for additional historical remarks
see \cite{p2}. Let us consider a regular pentagon whose edges have
length $1$ and let us consider the corresponding inscribed
pentagram, as in Figure~\ref{pentagramforkite}. It can be shown that
the length of the diagonal of the pentagon is equal to the {\em
golden ratio}, $\phi=\frac{1+\sqrt{5}}{2}=2\cos{\frac{\pi}{5}}$. The
polygon having vertices $A$, $B$, $E$ and $G$ is a Penrose kite, the
polygon having vertices $A$, $B$, $F$ and $G$ is a Penrose dart, and
their union is the Penrose thick rhombus having vertices $E$, $B$,
$F$ and $G$ (see Figure~\ref{kite}). Remark that the angles of the
kite measure $2\pi/5$ at the vertices $B$, $E$ and $G$ and $4\pi/5$
at the vertex $A$. Moreover, the longest edges, $EG$ and $EB$, and
the longest diagonal, $EA$, have the same length, which is $1$,
whilst the shortest edges  $AG$ and $AB$ have length $1/\phi$. The
angles of the dart measure $\pi/5$ at the vertices $G$ and $B$,
$2\pi/5$ at the vertex $F$ and $6/5 \pi$ at the vertex $A$.

\section{Quasilattices}

First of all let us recall the definition of quasilattice:
\begin{defn}[Quasilattice]{\rm
Let $V$ be a real vector space. A {\em quasilattice} in $V$ is the
span over $\Z$ of a set of $\R$--spanning vectors $V_1,\ldots,V_d$ of
$V$.}
\end{defn}
Notice that $\hbox{Span}_{\Z}\{V_1,\dots,V_d\}$ is a lattice if and
only if it admits a set of generators which is a basis of $V$.

\begin{figure}[h]
\begin{center}
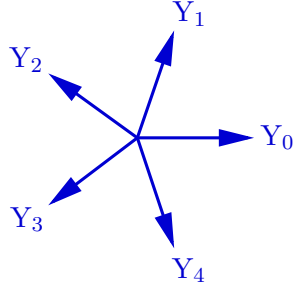
\end{center}
\caption{The star of vectors $Y_0,Y_1,Y_2,Y_3,Y_4$} \label{stargrid}
\end{figure}

It is easy to see that, in a suitably chosen coordinate system, the
unitary vectors
\begin{equation}\label{star}\left\{\begin{array}{l}
Y_1=(\cos{\frac{2\pi}{5}},\sin{\frac{2\pi}{5}})=\frac{1}{2}(\frac{1}{\phi},\sqrt{2+\phi})\\
Y_2=(\cos{\frac{4\pi}{5}},\sin{\frac{4\pi}{5}})=\frac{1}{2}(-\phi,\frac{1}{\phi}\sqrt{2+\phi})\\
Y_3=(\cos{\frac{6\pi}{5}},\sin{\frac{6\pi}{5}})=\frac{1}{2}(-\phi,-\frac{1}{\phi}\sqrt{2+\phi})\\
Y_4=(\cos{\frac{8\pi}{5}},\sin{\frac{8\pi}{5}})=\frac{1}{2}(\frac{1}{\phi},-\sqrt{2+\phi})
\end{array}
\right.\end{equation} are orthogonal to each of the four different
edges of the kite (cf. Figure~\ref{kite}). Now notice that
$Y_0=(1,0)$ is given by $Y_0=-(Y_1+Y_2+Y_3+Y_4)$ (see
Figure~\ref{stargrid}).
Therefore $\{Y_1,Y_2,Y_3,Y_4\}$ and $\{Y_0,Y_1,Y_2,Y_3,Y_4\}$
generate the same quasilattice, that we denote by $Q$, namely
\begin{equation}\label{qu}
Q=\hbox{Span}_{\Z}\{Y_0,Y_1,Y_2,Y_3,Y_4\}.\end{equation} The
quasilattice $Q$ is not a lattice, it is dense in $\R^2$ and a
minimal set of generators of $Q$ is made of four vectors.

The quasilattice $Q$ is naturally linked to the kite also in the
following sense. Consider the vectors dual to $Y_1,Y_2,Y_3,Y_4$;
they are given by
$$\left\{\begin{array}{l}
Y^*_1=\frac{1}{2}(-\sqrt{2+\phi},\frac{1}{\phi})\\
Y^*_2=\frac{1}{2}(-\frac{1}{\phi}\sqrt{2+\phi},-\phi)\\
Y^*_3=\frac{1}{2}(\frac{1}{\phi}\sqrt{2+\phi},-\phi)\\
Y^*_4=\frac{1}{2}(\sqrt{2+\phi},\frac{1}{\phi})
\end{array}
\right.$$
\begin{figure}[h]
\begin{center}
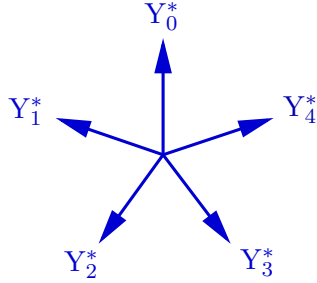
\end{center}
\caption{The star of vectors $Y^*_0,Y^*_1,Y^*_2,Y^*_3,Y^*_4$}
\label{stargriddual}
\end{figure}
Notice that, in the same coordinate system as above, the four edges
of the kite are parallel to these four vectors, and that its
vertices are contained in the quasilattice $R$ that they generate.
Notice that
$$Y^*_0=(0,1)$$ is given by
$$Y^*_0=-(Y^*_1+Y^*_2+Y^*_3+Y^*_4).$$ Therefore
$\{Y^*_1,Y^*_2,Y^*_3,Y^*_4\}$ and
$\{Y^*_0,Y^*_1,Y^*_2,Y^*_3,Y^*_4\}$ generate the same quasilattice.
We show the star of five vectors $Y^*_0,Y^*_1,Y^*_2,Y^*_3,Y^*_4$ in
Figure~\ref{stargriddual}.

\begin{figure}[h]
\begin{center}
\includegraphics{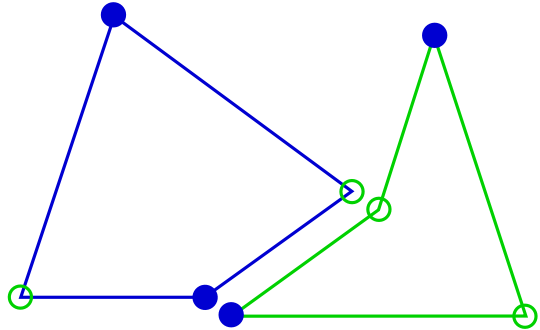}
\end{center}
\caption{Matching rules for the kite and dart tiling} \label{matchingruleskitedart}
\end{figure}

Let us now show that this connection between the quasilattice $R$
and the kite miraculously extends to any kite and dart tiling. We
recall that a kite and dart tiling is a tiling of the plane by kites
and darts that obey the matching rules shown in
Figure~\ref{matchingruleskitedart} (cf. \cite{pen1,pen2} and the
book by Senechal \cite{S} for a review on quasilattices and
nonperiodic tilings). There are uncountably many such tilings and
each of them is nonperiodic. Notice that the kite and dart can never
be joined to yield a thick rhombus, namely the configuration in
Figure~\ref{kite} is not allowed.
Consider now the vectors $Y_k^*$ and their opposites $-Y^*_k$, with
$k=0,\ldots,4$. For each vector $Y^*_k$ let $\D_k^+$ be the kite
such that, in the notations of Figure~\ref{kite}, $E$ coincides with
the origin and $A-E=Y^*_k$. We obtain in this way a star of five
kites. Analogously denote by $\D^-_k$ the kites corresponding to the
vectors $-Y^*_k$, thus obtaining a star of kites rotated by $\pi/5$
with respect to the first one.

Let us now consider any kite and dart tiling $\cal T$ with kites
having longest edge of length $1$. Denote by $AB$ one edge of the
tiling $\cal T$. From now on we will choose our coordinates so that
$A=O$ and so that $B-A$ is parallel to $Y^*_0$ with the same
orientation.

\begin{prop}\label{rotazioni}{\rm Let $\cal T$ be a kite and dart tiling
with kites having longest edge of length $1$. Then each kite of the
tiling is the translate of either a $\D_k^+$, $k=0,\ldots,4$ or a
$\D_k^-$, $k=0,\ldots,4$. Moreover each vertex of the tiling lies in
the quasilattice $R$.}\end{prop} \proof The argument is very simple.
Let $C$ be a vertex of the tiling that is different from $0$ and the
above vertex $B$. We can join $B$ to $C$ with a broken line made of
subsequent edges of the tiling. We denote the vertices of the broken
line thus obtained by $T_0=A,T_1=B,\dots,T_j,\dots,T_m=C$. The angle
of the broken line at each vertex $T_j$ is necessarily a multiple of
$\pi/5$ (see Section~\ref{figurekiteedart}). Therefore each vector
$V_j=T_j-T_{j-1}$ is either one of the vectors $\pm Y^*_k$,
$k=0,\ldots,4$, to account for the edges of length $1$, or one of
the vectors $\pm(Y_k^*+Y^*_{k+2})$, $k=0,\ldots,4$ (here
$Y^*_5=Y^*_0$ and $Y^*_6=Y^*_1$), to account for the edges of length
$1/\phi$. Since $C-A=T_m-T_0=V_m+\cdots+V_1$ our assertion is
proved: the vertex $C$ lies in $R$ and each kite having $C$ as
vertex is the translate of one of the ten kites $\D_k^+$ and
$\D_k^-$, $k=0,\ldots,4$. \qed
\begin{remark}\label{equivalenti}{\rm
Recall from \cite[Chapter~6, Section~1]{S} that a kite and dart
tiling gives rise to a rhombus tiling and viceversa. This can be
done by bisecting the tiles into isosceles triangles and then by
composing these triangles into thin and thick rhombuses. A
subdivided kite will become part either of a thick rhombus and a
thin rhombus, as shown in the local configuration (a) in
Figure~\ref{passaggio}, or of two thick rhombuses, as shown in the
local configuration (b). It can never become part of two thin
rhombuses, this is forbidden by the matching rules. All vertices of
a rhombus tiling are contained in the quasilattice $R$ (cf.
\cite[Proposition~1.3]{rhombus}).
\begin{figure}[h]
\begin{center}
\includegraphics{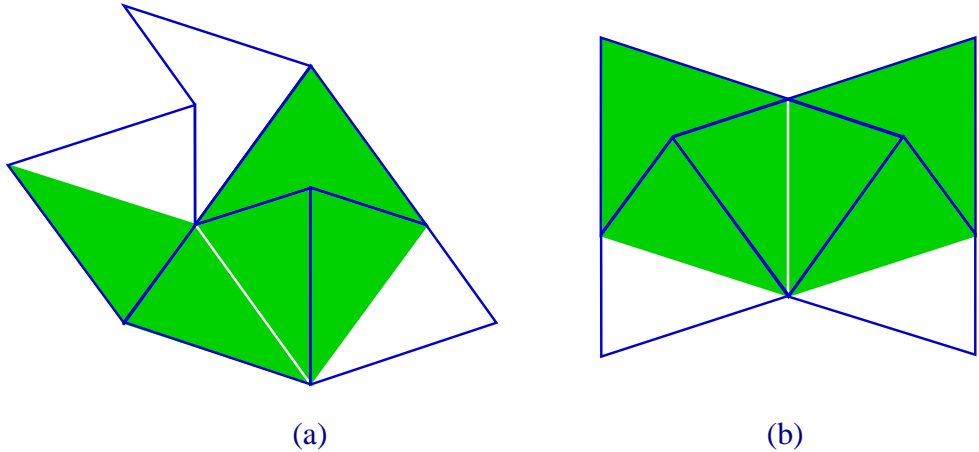}
\end{center}
\caption{Passing from a kite and dart tiling to a rhombus tiling}
\label{passaggio}
\end{figure}}
\end{remark}

\section{The Generalized Delzant Procedure}
We now outline the generalization of the Delzant procedure to
nonrational simple convex polytopes \cite{p1}. For the definition of
quasitori and their Hamiltonian actions we refer the reader to the
original article \cite{p1}, while for the definition of quasifolds
and related geometrical objects we refer to the Appendix.

Let us now recall what a simple convex polytope is.
\begin{defn}[Simple polytope] {\rm A dimension $n$ convex polytope $\D\subset\rndu$ is said to be {\em
simple} if there are exactly $n$ edges stemming from each
vertex.}\end{defn} Let us now consider a dimension $n$ convex
polytope $\D\subset\rndu$. If $d$ is the number of facets of $\D$,
then there exist elements $\xd$ in $\rn$ and $\ld$ in $\R$ such that
\begin{equation}\label{polydecomp}
\D=\bigcap_{j=1}^d\{\;\mu\in\rndu\;|\;\langle\mu,X_j\rangle\geq\lambda_j\;\}.
\end{equation}

\begin{defn}[Quasirational polytope]{\rm Let $Q$ be a quasilattice in
$\rn$. A convex polytope $\D\subset\rndu$ is said to be {\it
quasirational} with respect to $Q$ if the vectors $\xd$ can be
chosen in $Q$.}
\end{defn}
All polytopes in $\rndu$ are quasirational with respect to some
quasilattice $Q$; it is enough to consider the quasilattice that is
generated by the elements $\xd$ in (\ref{polydecomp}). Notice that
if the quasilattice is a honest lattice then the polytope is
rational.

In our situation we only need to consider the special case of simple
convex polytopes in $2$--dimensional space. Let $Q$ be a
quasilattice in $\rtwo$ and let $\D$ be a simple convex polytope in
the space $\rtwodu$ that is quasirational. Consider the space $\cd$
endowed with the standard symplectic form $\omega_0=\frac{1}{2\pi
i}\sum_{j=1}^d dz_j\wedge d\bar{z}_j$ and the standard action of the
torus $\td=\rd/\zd$:
$$
\begin{array}{cccccl}
\tau\,\colon& \td&\times&\cd&\longrightarrow& \cd\\
&((\et1,\ldots,\etd)&,&\vz)&\longmapsto&(\et1 z_1,\ldots, \etd z_d).
\end{array}
$$
This action is effective and Hamiltonian and its moment mapping is
given by
$$
\begin{array}{cccl}
J\,\colon&\cd&\longrightarrow &\rddu\\
&\vz&\longmapsto & \sum_{j=1}^d \zjs
e_j^*+\lambda,\quad\lambda\in\rddu \;\mbox{constant}.
\end{array}
$$
The mapping $J$ is proper and its image is the cone
$\cl=\lambda+\c0$, where $\c0$ denotes the positive orthant in the
space $\rddu$. Now consider the surjective linear mapping
\begin{eqnarray*}
\pi\,\colon &\rd \longrightarrow \rtwo,\\
&e_j \longmapsto X_j.
\end{eqnarray*}
Consider the dimension $2$ quasitorus $D=\rtwo/Q$. Then the linear
mapping $\pi$ induces a quasitorus epimorphism $\Pi\,\colon\,\td
\longrightarrow D$. Define now $N$ to be the kernel of the mapping
$\Pi$ and choose $\lambda=\sum_{j=1}^d \lambda_j e_j^*$. Denote by
$i$ the Lie algebra inclusion $\mbox{Lie}(N)\rightarrow\rd$ and
notice that $\Psi=i^*\circ J$ is a moment mapping for the induced
action of $N$ on $\cd$. Then the quasitorus $\td/N$ acts in a
Hamiltonian fashion on the compact symplectic quasifold
$M=\Psi^{-1}(0)/N$. If we identify the quasitori $D$ and $\td/N$
using the epimorphism $\Pi$, we get a Hamiltonian action of the
quasitorus $D$ whose moment mapping has image equal to
${(\pi^*)}^{-1}(\cl\cap\ker{i^*})=
{(\pi^*)}^{-1}(\cl\cap\mbox{im}\,\pi^*)= {(\pi^*)}^{-1}(\pi^*(\D))$
which is exactly $\D$. This action is effective since the level set
$\Psi^{-1}(0)$ contains points of the form $\vz\in\cd$, $z_j\neq0$,
$j=1,\ldots,d$, where the $T^d$-action is free. Notice finally that
$\dim{M}=2d-2\dim{N}= 2d-2(d-2)=4=2\dim{D}$.

\begin{remark}{\rm This construction depends on two arbitrary
choices: the choice of the quasilattice $Q$ with respect to which
the polytope is quasirational, and the choice of the
inward--pointing vectors $\xd$ in $Q$.}
\end{remark}
\begin{remark}\label{connesso}{\rm
It is easy to show that if the vectors $\xd$ are generators of the
quasilattice $Q$, then $N=\exp(\n)$ and is therefore connected.}
\end{remark}

\section{The Kite from a Symplectic Viewpoint}\label{kitesimplettico}
\begin{figure}[h]
\begin{center}
\includegraphics{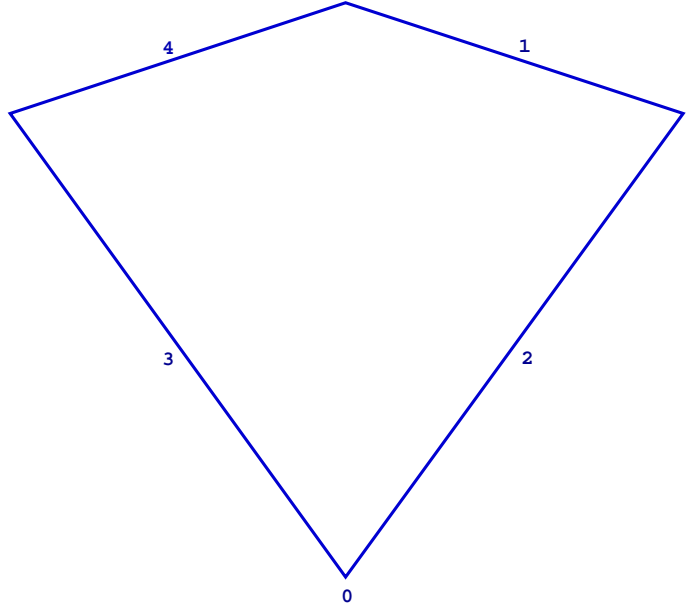}
\end{center}
\caption{The kite} \label{kite2}
\end{figure}

Let us apply the generalized Delzant procedure to the kite $\D_0^+$.
Let us label its edges with the numbers $1,2,3,4$, as in
Figure~\ref{kite2}.
Our choice of inward--pointing vectors is given by $X_1=-Y_1$,
$X_2=Y_2$, $X_3=-Y_3$, $X_4=Y_4$. Remark that the vectors
$X_1,X_2,X_3,X_4$ generate the quasilattice $Q$ defined by
(\ref{qu}), and that the kite is quasirational with respect to $Q$.
It is easy to see that the constants $\lambda_1,\ldots,\lambda)4$ in
\ref{polydecomp} are given by
 $\lambda_1=\lambda_4=
-\frac{1}{2}\sqrt{2+\phi}$, and $\lambda_2=\lambda_3=0$. Let us
consider the linear mapping defined by
$$
\begin{array}{cccc}
\pi\,\colon\,&\R^4&\rightarrow&\R^2\\
&e_i&\mapsto &X_i.
\end{array}
$$
Consider now the subgroup $N=\{\,\exp(X)\in T^4\,|\,X\in\R^4 ,
\pi(X)\in Q\,\}$. It is easy to see, using the relations
$$Y_2=-Y_1-\phi Y_4$$
$$Y_3=-\phi Y_1-Y_4,$$ that
its Lie algebra is given by
$$\n=\left\{\,X\in\R^4\,|\, X=\left(-s+\phi t,s,t,-t+\phi s\right),
s,t \in\R \,\right\}.$$ Therefore, by Remark~\ref{connesso} we have
that
$$N=\exp(\n)=\left\{\,\exp(X)\in T^4\,|\, X=\left(-s+\phi t,s,t,-t+\phi s\right),
s,t \in\R\,\right\}.$$ We will be needing the following bases for
$\n$
$$
\begin{array}{ll}
B_{12}&=\{(1,1/\phi,1,0),(-1/\phi,1/\phi,0,1)\}\\
B_{34}&=\{(1,0,1/\phi,-1/\phi),(0,1,1/\phi,1)\}
\end{array}
$$
and the following identity for the golden ratio
$$\phi=1+\frac{1}{\phi}\;.$$ Let us consider $\psi$, the moment mapping of
the induced $N$--action, and let us write it down in two different
ways, relatively to the two different bases above:
\begin{eqnarray}
\psi(z_1,z_2,z_3,z_4)&=&\left(|z_1|^2+\frac{1}{\phi}|z_2|^2+|z_3|^2-\phi\sigma, -\frac{1}{\phi}|z_1|^2+\frac{1}{\phi}
|z_2|^2+|z_4|^2-\frac{\sigma}{\phi}\right)\label{12}\\
&=&\left(|z_1|^2+\frac{1}{\phi}|z_3|^2-\frac{1}{\phi}|z_4|^2-\frac{\sigma}{\phi},
|z_2|^2+\frac{1}{\phi}|z_3|^2+|z_4|^2-\sigma\phi\label{34}\right)
\end{eqnarray}
where we write $\sigma=\frac{1}{2\phi}\sqrt{2+\phi}$ for brevity.

According to the generalized Delzant procedure, the quotient space
$M^+_0=\Psi^{-1}(0)/N$ is the quasifold associated to the kite
$\D_0^+$. The quasitorus $D^2=\R^2/Q$ acts on $M_0^+$ in a
Hamiltonian fashion, with image of the corresponding moment mapping
yielding exactly the kite $\D_0^+$.

The atlas that defines the quasifold structure of $M_0^+$ is given
by four charts, each of which corresponds to a vertex of the kite.
We will describe just two of them.

First consider the vertex that is the intersection of the edges
labeled $1$ and $2$. Consider the open subset of $\C^2$ given by
$$\tilde{U}_{12}=\left\{\,(z_1,z_2)\in\C^2\,|\,|z_1|^2+\frac{1}{\phi}|z_2|^2<\phi\sigma,\,-|z_1|^2+|z_2|^2<\sigma\,\right\}.$$
We now use (\ref{12}) to construct the following slice of
$\Psi^{-1}(0)$ that is transversal to the $N$--orbits
$$\begin{array}{ccc}
\tilde{U}_{12}& \stackrel{\tilde{\tau}_{12}}{\longrightarrow}&
\{\vz\in\Psi^{-1}(0)\;|\;z_3\neq0,z_4\neq0\}\\
(z_1,z_2)&\longmapsto&\left(z_1,z_2,\sqrt{\phi\sigma-|z_1|^2-\frac{1}{\phi}|z_2|^2},\sqrt{\frac{1}{\phi}(\sigma+|z_1|^2-|z_2|^2})\right)
\end{array}
$$
which induces the homeomorphism
$$
\begin{array}{ccc}
\tilde{U}_{12}/\Gamma_{12}&\stackrel{\tau_{12}}{\longrightarrow}& U_{12}\\
\,[(z_1,z_2)]&\longmapsto&[{\tilde{\tau}}_{12}(z_1,z_2)]
\end{array},
$$
where the open subset $U_{12}$ of $M^+_0$ is the quotient
$$\{\vz\in\Psi^{-1}(0)\;|\;z_3\neq0,z_4\neq0\}/N$$
and the discrete group $\Gamma_{12}$ is given by
\begin{equation}\label{gamma12}
\Gamma_{12}=\left\{\,(e^{-2\pi i\frac{1}{\phi}h},e^{2\pi i \frac{1}{\phi}(h+k)}\in T^2\;|\; h,k\in\Z\right\}.
\end{equation}
The triple $(U_{12},\tau_{12},\tilde{U}_{12}/\Gamma_{12})$ defines a
chart for $M^+_0$.

To construct a second chart we consider the vertex that is given by
the intersection of the edges $3$ and $4$. Consider the open subset
of $\C^2$ given by
$$\tilde{U}_{34}=\left\{\,(z_3,z_4)\in\C^2\,|\,|z_3|^2-|z_4|^2<{\sigma},
\frac{1}{\phi}|z_3|^2+|z_4|^2<\sigma\phi\,\right\}.$$
We now use (\ref{34}) to construct the following slice
$$\begin{array}{ccc}
\tilde{U}_{34}& \stackrel{\tilde{\tau}_{34}}{\longrightarrow}&
\{\vz\in\Psi^{-1}(0)\;|\;z_1\neq0,z_2\neq0\}\\
(z_3,z_4)&\longmapsto&\left(\sqrt{\frac{1}{\phi}\left(\sigma-|z_3|^2+|z_4|^2\right)},\sqrt{\sigma\phi-\frac{1}{\phi}|z_3|^2-|z_4|^2},z_3,z_4\right)
\end{array}
$$
which induces the homeomorphism
$$
\begin{array}{ccc}
\tilde{U}_{34}/\Gamma_{34}&\stackrel{\tau_{34}}{\longrightarrow}& U_{34}\\
\,[(z_3,z_4)]&\longmapsto&[{\tilde{\tau}}_{34}(z_3,z_4)]
\end{array},
$$
where the open subset $U_{34}$ of $M^+_0$ is the quotient
$$\{\vz\in\Psi^{-1}(0)\;|\;z_1\neq0,z_2\neq0\}/N$$
and the discrete group $\Gamma_{23}$ is given by
$$\Gamma_{34}=\left\{\,(e^{2\pi i \frac{1}{\phi}(h+k)},e^{-2\pi i \frac{1}{\phi}h})\in T^2\;|\; h,k\in\Z\right\}.$$
This yields our second chart
$(U_{34},\tau_{34},\tilde{U}_{34}/\Gamma_{34})$. The two other
charts are constructed in a similar way.

In order to show that the four charts are compatible we need to show
that the changes of charts are well defined for each pair of
overlapping charts. Let us see this in detail for the pair of charts
$U_{12}$ and $U_{34}$. We prove that for each $m\in U_{12}\cap
U_{34}$ we can take as subset $U^m_{12\,34}$ of
Definition~\ref{compatiblecharts} the connected open subset
$U_{12}\cap U_{34}$. Observe first that
$$\tau^{-1}_{12}(U_{12}\cap U_{34})=
\{(z_1,z_2)\in\tilde{U}_{12}\;|\;z_1\neq0,z_2\neq0\}/\G_{12}
$$
and that, in the same way,
$$\tau^{-1}_{34}(U_{12}\cap U_{34})=
\{(z_3,z_4)\in\tilde{U}_{34}\;|\;z_3\neq0,z_4\neq0\}/\G_{34}.
$$
Let us introduce
$$\tilde{V}_{12}=\{(z_1,z_2)\in\tilde{U}_{12}\;|\;z_1\neq0,z_2\neq0\}$$
and
$$\tilde{V}_{34}=\{(z_3,z_4)\in\tilde{U}_{34}\;|\;z_3\neq0,z_4\neq0\}.$$
We need to prove that the mapping
$$\tau^{-1}_{34}\circ\tau_{12}\,\colon\,\tilde{V}_{12}/\G_{12}\longrightarrow
\tilde{V}_{34}/\G_{34}$$ is a diffeomorphism of the universal
covering models of the given models $\tilde{V}_{12}/\G_{12}$ and
$\tilde{V}_{34}/\G_{34}$. This means, by definition of model
diffeomorphism, that $\tau^{-1}_{34}\circ\tau_{12}$  lifts to a
diffeomorphism between the appropriate open subsets of $\R^4$. The
set $\tilde{V}_{12}$ has first fundamental group isomorphic to
$\Z\times\Z$ and its universal covering is the open subset
$$\vsh_{12}=\left\{(\rho_1,\theta_1,\rho_2,\theta_2)\in (\R^+\times\R)^2\;|\;
|\rho_1|^2+\frac{1}{\phi}|\rho_2|^2<\phi\sigma,\,-|\rho_1|^2+|\rho_2|^2<\sigma\,\right\}.$$
The discrete group $\Gsh_{12}\simeq\Z^4$
acts on $\vsh_{12}$ in the following way:
\begin{equation}\label{azionedisharp12}
\xymatrix@R=1pc{
\Gsh_{12}\times\vsh_{12}\ar[r]&\vsh_{12}\\
(m,n,h,k)\times(\rho_1,\theta_1,\rho_2,\theta_2)\ar@{|->}[r]&
(\rho_1,\theta_1-\frac{1}{\phi} h+m,\rho_2,\theta_2+\frac{1}{\phi}(h+k)+n).
}
\end{equation}
This action on $\vsh_{12}$ satisfies the hypotheses of
Definition~\ref{model}. The group $\Gsh_{12}$ is an extension of the
group $\G_{12}$ by the group $\hbox{Aut}(\vsh_{12}|\vt_{12})$
$$1\longrightarrow \hbox{Aut}(\vsh_{12}|\vt_{12})\longrightarrow\Gsh_{12}
\longrightarrow\G_{12}\longrightarrow 1.$$ Therefore
$\vsh_{12}/\Gsh_{12}$ is the universal covering model of
$\tilde{V}_{12}/\G_{12}$. Analogously we consider the universal
covering  model of $\tilde{V}_{34}/\G_{34}$ given by the quotient of
the open set
$$\vsh_{34}=\left\{(\rho_3,\theta_3,\rho_4,\theta_4)\in (\R^+\times\R)^2\;|\;
\rho_3^2-|\rho_4|^2<{\sigma},\frac{1}{\phi}\rho_3^2+|\rho_4|^2<\sigma\phi
\right\}$$
by the following action of the group $\Gsh_{34}\simeq\Z^4$:
\begin{equation}\label{azionedisharp34}
\xymatrix@R=1pc{
\Gsh_{34}\times\vsh_{34}\ar[r]&\vsh_{34}\\
(m,n,h,k)\times(\rho_3,\theta_3,\rho_4,\theta_4)\ar@{|->}[r]&
(\rho_3,\theta_3+\frac{1}{\phi}(h+k)+m,\rho_4,\theta_4-\frac{1}{\phi}h+n).
}
\end{equation}
The mapping $(\tau^{-1}_{34}\circ\tau_{12})^{\sharp}$
$$
\xymatrix@R=1pc{
\vsh_{12}\ar[r]&\vsh_{34}\\
(\rho_1,\theta_1,\rho_2,\theta_2)\ar@{|->}[r]&
\left(\sqrt{\phi\sigma-\rho_1^2-\frac{1}{\phi}\rho_2^2},-\frac{1}{\phi}(\theta_1+\theta_2),
\sqrt{\frac{1}{\phi}(\sigma+\rho_1^2-\rho_2^2)},\frac{1}{\phi}\theta_1-\theta_2\right)
}
$$
is a diffeomorphism. In particular, it is a lift of
$\tau^{-1}_{34}\circ\tau_{12}$, namely the following diagram is
commutative:
\begin{equation}\label{changeofcharts1234}
\xymatrix@C=4pc{
\vsh_{12}\ar[d]\ar[r]^{(\tau^{-1}_{34}\circ\tau_{12})^{\sharp}}
&\vsh_{34}\ar[d]\\
\vsh_{12}/\Gsh_{12}\ar[r]^{\tau^{-1}_{34}\circ\tau_{12}}&\vsh_{34}/\Gsh_{34}}
\end{equation}
This shows that the two charts are compatible.
\begin{remark}\label{impossibile}{\rm
The group isomorphism $F\,\colon\,\Gsh_{12}\rightarrow\Gsh_{34}$
that makes $(\tau^{-1}_{34}\circ\tau_{12})^{\sharp}$ equivariant is
given by $F(m,n,h,k)=(-k,-n-h,k-m,-n)$. This implies that
$(\tau^{-1}_{34}\circ\tau_{12})^{\sharp}$ is not equivariant with
respect to the actions of $\hbox{Aut}(\vsh_{12}|\vt_{12})$ and
$\hbox{Aut}(\vsh_{34}|\vt_{34})$. Therefore there does not exist a
mapping $\widetilde{(\tau^{-1}_{34}\circ\tau_{12})}$ that makes the
below diagram commutative
$$
\xymatrix@C=4pc{
\vsh_{12}\ar[d]\ar[r]^{(\tau^{-1}_{34}\circ\tau_{12})^{\sharp}}
&\vsh_{34}\ar[d]\\
\vt_{12}\ar[r]^{\widetilde{(\tau^{-1}_{34}\circ\tau_{12})}}&\vt_{34}
}
$$
}
\end{remark}
\medskip
\noindent We proceed in the same way for the other pairs of
overlapping charts. The four charts turn out to be compatible, thus
defining on $M^+_0$ a quasifold structure.

We now describe explicitly the symplectic structure on $M^+_0$
induced by the reduction procedure. To define a symplectic structure
on a quasifold we first need to define a symplectic structure on
each chart, and then require that the different structures behave
well under changes of charts. Consider, for example, the chart
$U_{12}\simeq \tilde{U}_{12}/\G_{12}.$ A symplectic form here is
given by a symplectic form on $ \tilde{U}_{12}$ which is invariant
under the action of $\G_{12}$. We take the restriction  of the
standard symplectic form on $\C^2$ to $ \tilde{U}_{12}$ and we
denote it by $\tilde{\omega}_{12}$. We do the same for the three
other charts. Consider now the changes of charts, for example the
one described above, and let $\omega^{\sharp}_{12}$ and
$\omega^{\sharp}_{34}$ be the pullbacks of the forms
$\tilde{\omega}_{12}$ and $\tilde{\omega}_{34}$ to $\vsh_{12}$ and
$\vsh_{34}$ respectively. It is easy to check that the pullback of
$\omega^{\sharp}_{34}$ via the mapping
$(\tau^{-1}_{34}\circ\tau_{12})^{\sharp}$ is exactly the form
$\omega^{\sharp}_{12}$. The same argument applies to the other
changes of charts. This symplectic structure is the one that is
induced by the reduction procedure, namely the pullback of $\omega$
via the projection to the quotient coincides with the pullback of
the standard form on $\C^4$ via the inclusion mapping:
$$
\xymatrix{
\Psi^{-1}(0)\ar@{^{(}->}[r]\ar[d]&\C^4\\
\Psi^{-1}(0)/N
}
$$

\section{Global Symplectic Interpretation of a Kite and Dart Tiling}
Recall that we denoted by $M_0^+$ the symplectic quasifold
associated to the kite $\D_0^+$. Consider the ten distinguished
kites $\D_k^+$ and $\D_k^-$, $k=0,\dots,4$. Notice that each of
these kites has a natural choice of inward--pointing orthogonal
vectors; these are $-Y_{k+1}$, $Y_{k+2}$, $-Y_{k+3}$, $Y_{k+4}$ for
$\D_k^+$, and $Y_{k+1}$, $-Y_{k+2}$, $Y_{k+3}$, $-Y_{k+4}$ for
$\D_k^-$. Consider now a kite and dart tiling with longest edges of
length $1$. Remark that, by Proposition~\ref{rotazioni}, in our
choice of coordinates, {\em each} of its kites can be obtained by
translation from one of the $10$ kites $\D_k^+$ and $\D_k^-$. We can
then prove the following
\begin{thm}\label{uguali}{\rm The compact connected symplectic quasifold
corresponding to each kite of a kite and dart tiling with longest
edges of length $1$ is given by $M_0^+$.}\end{thm} \proof Observe
that, for each $k=1,\ldots,4$, there exists an orthogonal
transformation $P$ of $\R^2$ leaving the quasilattice $Q$ invariant,
that sends the orthogonal vectors relative to the kite $\D_k^+$ to
the orthogonal vectors relative to the kite $\D_0^+$, and such that
the dual transformation $P^*$ sends the kite $\D^+_0$ to the kite
$\D_k^+$. The same is true for the kites $\D^-_k$, with
$k=1,\ldots,4$. This implies that the reduced space corresponding to
each of the $10$ kites $\D_k^+$ and $\D_k^-$, with the choice of
orthogonal vectors and quasilattice specified above, is exactly
$M^+_0$. This yields a unique symplectic quasifold, $M^+_0$, for all
the kites considered. Finally, it is straightforward to check that
translating the kites $\D_k^+$ and $\D_k^+$ does not produce any
change in the corresponding quotient spaces, therefore, by
Proposition~\ref{rotazioni} we are done. \qed
\begin{remark}{\rm We recall that, with our choices of quasilattice and inward pointing vectors,
we obtain a unique symplectic quasifold $M_0^+$ for any tiling
having the same edge lengths. We suggest a modification of our
approach that allows one to distinguish between different tilings.
In \cite{dbr} De Bruijn gives a construction that associates to each
suitably regular point of the hyperplane
$H=\left\{(\gamma_0,\ldots,\gamma_4)\in\R^5\;|\;\sum_{j=0}^4\gamma_j=0\right\}$
a unique Penrose rhombus tiling and therefore a unique kite and dart
tiling. Let us consider a regular point $(\gamma_0,\ldots,\gamma_4)$
in $H$. By \cite[Theorem~9.2]{dbr} it is always possible to assume,
up to translation of the corresponding tiling, that $\gamma_j>0$ for
$j=1,2,3,4$ and $\gamma_0<0$. We can therefore apply the generalized
Delzant procedure to any kite of the corresponding tiling with
respect to the quasilattice that is generated by $\{-\gamma_0 Y_0,
\gamma_1 Y_1,\gamma_2 Y_2,\gamma_3 Y_3,\gamma_4 Y_4\}$. With this
choice the corresponding symplectic quasifold does keep track of the
quintuple $(\gamma_0,\ldots,\gamma_4)$ that characterizes the kite
and dart tiling. We plan to investigate this approach in future
work. }
\end{remark}

\section{The quasifold $M^+_0$ is not a global quotient}

The present section is devoted to the proof of the following theorem:

\begin{thm}\label{noneglobale}{\rm The quasifold $M_0^+$ is not a global quotient.}\end{thm}
For the proof we need some intermediate results that are better
proved in a more general and simple context. These results are
summarized in the following Lemmas and Remarks:

\begin{lemma}\label{lemma1}{\rm Let $(\ut/\G, p, \ut)$ and $(\vt/\D,q,\vt)$ be  two models and
let $f$ be a diffeomorphism of their universal covering models.
Suppose that $\ut$ is simply connected and that there exists a point
$\tilde{u}_0\in\ut$ such that the isotropy of $\G$ at $\tilde{u}_0$,
$\G_{\tilde{u}_0}$, is the whole group $\G$. Then $\vt$ is itself
simply connected, $f$ is a diffeomorphism between the two given
models and $\G$ and $\D$ are isomorphic.}
\end{lemma}
\proof By hypothesis we have the following diagram:
$$
\xymatrix{
&\vsh\ar[d]^{\pi}\\
\ut\ar[ur]^{\fsh}\ar[d]^{p}&\vt\ar[d]^{q}\\
\ut/\G\ar[r]^{f}&\vsh/\Dsh
}
$$
By Lemma~\ref{green} we have an isomorphism
$F\,\colon\,\G\rightarrow\D^{\sharp}$ that makes the lift $\fsh$
equivariant. Denote by $v^{\sharp}_0$ the point $\fsh(\tilde{u}_0)$.
Take $\gamma\in\G$; we have, by assumption,
$\gamma\tilde{u}_0=\tilde{u}_0$, therefore
$\fsh(\tilde{u}_0)=\fsh(\gamma\tilde{u}_0)=F(\gamma)\fsh(\tilde{u}_0)$.
Hence $F\,\colon\,\G\rightarrow\D^{\sharp}_{v^{\sharp}_0}$ and
$\D^{\sharp}=\D_{v^{\sharp}_0}^{\sharp}$. We therefore have the
following commutative diagram:
$$
\xymatrix{
&&1&&\\
1\ar[r]&\hbox{Aut}(\vsh|\vt)\ar[r]&\D^{\sharp}\ar[r]\ar[u]&\D\ar[r]&1\\
&1\ar[r]&\D_{v^{\sharp}_0}^{\sharp}\ar[r]\ar[u]&\D_{\pi({v^{\sharp}_0})}\ar[r]\ar[u]&1\\
&&1\ar[u]&1\ar[u]&
}
$$
which implies that $\hbox{Aut}(\vsh|\vt)$ is trivial. Therefore
$\vt$ is itself symply connected and $f$ is a diffeomorphism between
the two given models. This also implies that $\G$ and $\D$ are
isomorphic.\qed
\begin{lemma}\label{lemma3}{\rm
Let $(\ut/\G, p, \ut)$ and $(\vt/\D,q,\vt)$ be  two models, and let
$(\ud/\Gd,p\circ\pi, \ud)$ and $(\vd/\Dd,q\circ\rho,\vd)$ be two
covering models (here $\ud$ and $\vd$ are not necessarily simply
connected). Let $f\,\colon\,\ut/\G\longrightarrow\vt/\D$ be a
homeomorphism such that there are a diffeomorphism
$\ft_1\,\colon\,\ut\rightarrow\vt$ and a diffeomorphism
$\fd_2\,\colon\,\ud\rightarrow\vd$ that are both lifts of $f$. Then
$\rho\circ\fd_2$ descends to a diffeomorphism
$\ft_2\,\colon\,\ut\rightarrow\vt$. Moreover, there exists
$\delta\in\D$ such that $\ft_2=\delta\cdot\ft_1$.}\end{lemma} \proof
Consider the diagram:
$$
\xymatrix{
\ud\ar[d]^{\pi}\ar[r]^{\fd_2}\ar@/^2pc/[r]^{\fd_1}&\vd\ar[d]^{\rho}\\
\ut\ar[r]^{\ft_1}\ar[d]^{p}&\vt\ar[d]^{q}\\
\ut/\G\ar[r]^{f}&\vt/\D
}
$$
where $\fd_1$ is a lift of $\ft\circ\pi$ to $\ud$; such a lift
exists since $(\ft_1\circ\pi_u)_*(\pi_1(\ud))\subset \pi_1(\vt)$ and
is a diffeomorphism. By Lemma~\ref{green} there exists
$\delta^{\dagger}\in\D^{\dagger}$ such that
$\fd_1=\delta^{\dagger}\cdot\fd_2$. Therefore $\rho\circ\fd_2$
descends to a diffeomorphism $\ft_2\,\colon\,\ut\rightarrow\vt$ and
there exists $\delta\in\D$ such that $\ft_1=\delta\cdot\ft_2$. .\qed
\begin{remark}\label{remark2}{\rm
Let $(\ut/\G, p, \ut)$ and $(\vt/\D,q,\vt)$ be  two models, let
$(\ush/\Gsh,p\circ\pi,\ush)$ and $(\vsh/\Dsh,q\circ\rho,\vsh)$ be
their universal covering models, and let
$f\,\colon\,\ush/\Gsh\longrightarrow\vsh/\Dsh$ be a diffeomorphism.
Let $\fsh\,\colon\,\ush\rightarrow\vsh$ be the lift of $f$. Now take
an open connected subset $U_1\subset\ut/\G$. Let
$(\ut_1/\G_1,p_1,\ut_1)$ be an induced model. Denote by $i$ the
inclusion $\ut_1\hookrightarrow\ut$, suppose that $\G_1=\G$ and that
$i_*\,\colon\,\pi_1(\ut_1)\rightarrow\pi_1(\ut)$ is an isomorphism.
Then $\pi^{-1}(\ut_1)$ is the universal covering of $\ut_1$ and
$\pi^{-1}(\ut_1)/\Gsh$ is the universal covering model of
$\ut_1/\G$.}
\end{remark}
\begin{lemma}\label{lemma2}{\rm Let $(\ut/\G, p, \ut)$ and $(\vt/\D,q,\vt)$ be  two models,
let $(\ush/\Gsh,p\circ\pi,\ush)$ and $(\vsh/\Dsh,q\circ\rho,\vsh)$
be their universal coverings models and let
$f\,\colon\,\ush/\Gsh\longrightarrow\vsh/\Dsh$ be a diffeomorphism.
Let $\ut_1/\G_1\subset\ut/\G$ be as in Remark~\ref{remark2} and let
$\fsh_1$ be the restriction of $\fsh$ to $\pi^{-1}(\ut_1)$. If
$\rho\circ\fsh_1$ descends to a diffeomorphism
$\ft_1\,\colon\,\ut_1\rightarrow(\rho\circ\fsh_1)(\pi^{-1}(\ut_1))$,
then $\rho\circ\fsh$ descends to a diffeomorphism
$\ft\,\colon\,\ut\rightarrow\vt$. }\end{lemma} \proof Let
$F\,\colon\,\Gsh\rightarrow\Dsh$ be the group isomorphism that makes
$\fsh$ equivariant. Then the group isomorphism $F_1$ that makes
$\fsh_1$ equivariant is the restriction of $F$ to $\Gsh_1$; however,
$\Gsh_1=\Gsh$, therefore $F=F_1$. By assumption,
$F_1\,\colon\,\hbox{Aut}(\pi^{-1}(\ut_1)|\ut)\rightarrow
\hbox{Aut}(\vsh_1|\vt_1)$ is an isomorphism. On the other hand, we
can deduce from
$\hbox{Aut}(\pi^{-1}(\ut_1)|\ut)=\hbox{Aut}(\ush|\ut)$ that
$F\,\colon\,\hbox{Aut}(\ush,\ut)\rightarrow\hbox{Aut}(\vsh,\vt)$ is
an isomorphism. The existence of $\ft$ then follows. \qed

We are now ready to prove the main result:\\
{\mbox{\bf Proof of Theorem~\ref{noneglobale}}. We need to prove
that there does not exist a quasifold diffeomorphism $f$ from
 $M_0^+$ onto a global quotient $N=\tilde{N}/\D$.
Suppose that such a  diffeomorphism $f$ exists and let us show that
this leads to a contradiction. Let $m^1\in U_{12}\cap U_{34}$ and
let $n^1=f(m^1)$. First of all observe that, if $W$ is an open
subset of $M_0^+$ containing $m^1$ and such that
$\wt_{12}=(\tau_{12}\circ p_{34})^{-1}(W)$ is a product of two open
annuli, then $W\subset U_{12}\cap U_{34}$, and
$$\wt_{34}=(\tau_{12}\circ p_{34})^{-1}(W)$$ is itself a product of two open annuli.
Moreover, $\wsh_{12}=\pi_{12}^{-1}(\wt_{12})$ is the universal
covering of $\wt_{12}$ and so is $\wsh_{34}=\pi_{34}^{-1}(\wt_{34})$
with respect to $\wt_{34}$. Therefore the change of charts
$g=\tau_{34}\circ\tau_{12}^{-1}$ lifts to a diffeomorphism
$g^{\sharp}_W$ from $\whs_{12}$ onto $\wsh_{34}$, which is simply
the restriction to $\whs_{12}$ of the change of chart
$g^{\sharp}\,\colon\,\vsh_{12}\rightarrow\vsh_{34}$. It is crucial
to notice that, by Remark~\ref{impossibile},
 $\pi_{34}\circ g^{\sharp}_W$ {\em does not} descend to a diffeomorphism $\tilde{g}$ from $\wt_{12}$ to $\wt_{34}$.
Our strategy to prove that the diffeomorphism $f$ cannot exist is to
show that, if it does, then $\tilde{g}$ exists, leading to a
contradiction. Consider first the chart $U_{12}$, take a point
$(z^1_1,z^1_2)\in\vt_{12}$ such that $(\tau_{12}\circ
p_{12})((z^1_1,z^1_2))=m^1$ and denote by $m^0$ the point
$(\tau_{12}\circ p_{12})(0)$. The curve $t(z^1_1,z^1_2)$, with
$t\in[0,1]$, projects to $m^t=(\tau_{12}\circ
p_{12})(t(z^1_1,z^1_2))$. Notice that we can deduce from
(\ref{gamma12}) that any $\G_{12}$-invariant open subset of $V_{12}$
containing the point $(z^1_1,z^1_2)$ must contain the product of
circles
$\{(z_1,z_2)\in\vt_{12}\;|\;|z_1|=|z^1_1|,\,|z_2|=|z^2_1|\}$. Let us
now take our diffeomorphism $f\,\colon\, M_0^+\rightarrow N$. By
Definition~\ref{diffeo}, Remark~\ref{tuttelecarte} and
Proposition~\ref{restrizionedimappe}, for each $t\in(0,1]$ we can
choose a connected open subset $W$ around $m^t$ such that
$(\tau_{12}\circ p_{34})^{-1}(W)$ is a product of two open annuli
and $f\circ\tau_{12}\,\colon\,\tau_{12}^{-1}(W)\rightarrow f(W)$ is
a diffeomorphism of the universal covering models of the induced
models. When $t=0$ the open set $W$ can be chosen in such a way that
$(\tau_{12}\circ p_{34})^{-1}(W)$ is a product of two open disks. We
can cover the curve $m_t$ by a finite number of such $W_t$'s. Let us
denote this covering by $W_j$, $j=0,\ldots,k$, where $W_0$ is a
product of open disks and $W_j\cap W_{j+1}\neq0$.

It is convenient now to divide the last part of the proof in successive steps:\\
\underline{Step 1}: We start by considering $W_0$. Notice that the
isotropy of $\G_{12}$ at $0$ is exactly $\G_{12}$. Therefore, by
Lemma~\ref{lemma1}, we can conclude that the homeomorphism
$h_0=f\circ\tau_{12}$, defined on $\tau_{12}^{-1}(W_0)$, is a
diffeomorphism from the model  $(\tau_{12}\circ
p_{12})^{-1}(W_0)/\G_{12}$ onto a model
induced by $f(W_0)\subset\tilde{N}/\D$.\\
\underline{Step 2}: Consider the homeomorphism $h_1=f\circ\tau_{12}$
defined on $\tau_{12}^{-1}(W_1)$. By construction  $h_1$ is a
diffeomorphims of the universal covering models of the induced
models. We therefore have a diagram of the kind:
$$\xymatrix{
\wsh_{1}\ar[r]^{h^{\sharp}_1}\ar[d]^{\pi_1}&\vsh_1\ar[d]_{\rho_1}\\
\wt_{1}\ar[d]^{p_1}&\vt_1\ar[d]^{q_1}\\
\tau_{12}^{-1}(W_1)\ar[r]^{h_1}&f(W_1)
}
$$
Consider the restriction of $h_1$ to $\tau_{12}^{-1}(W_0\cap W_1)$.
This restriction admits a lift, defined on $(\pi_1\circ
p_{1})^{-1}(\tau_{12}^{-1}(W_0\cap W_1))$, which is simply the
restriction of $h^{\sharp}_1$. By Remark~\ref{tuttistessopi} all of
the models induced by $f(W_0\cap W_1)\subset \tilde{N}/\D$ are
diffeomorphic. Then, by  Step 1, the restriction of $h_1$ admits
another lift; this one is the restriction to
$p_{1}^{-1}(\tau_{12}^{-1}(W_0\cap W_1))$ of the lift of $h_0$.
Therefore, by Lemma~\ref{lemma3}, the restriction of $\rho_1\circ
h_1^{\sharp}$ to $(\pi_1\circ p_{1})^{-1}(\tau_{12}^{-1}(W_0\cap
W_1))$
descends to a diffeomorphism defined on $p_{1}^{-1}(\tau_{12}^{-1}(W_0\cap W_1))$.\\
\underline{Step 3}: Consider $W_0\cap W_1\subset W_1$. We are now in
the position of applying Lemma~\ref{lemma2} to the homeomorphism
$f\circ\tau_{12}$ defined on $\tau_{12}^{-1}(W_1)$, which is
therefore a diffeomorphism of the model $(\tau_{12}\circ
p_{12})^{-1}(W_1)/\G_{12}$
to a  model induced by $f(W_1)\subset\tilde{N}/\G_{\alpha_1}$.\\
\underline{Step 4}: By applying Step 3 to the other successive
intersections we show that $f\circ\tau_{12}$ is a diffeomorphism of
the model $(\tau_{12}\circ p_{12})^{-1}(W_k)/\G_{12}$
onto a model induced by $f(W_k)\subset\tilde{N}/\D$.\\
\underline{Step 5}: By applying Steps 1 to 4 to the chart $U_{34}$
we show that $f\circ\tau_{34}$ is a diffeomorphism of the model
$(\tau_{34}\circ p_{34})^{-1}(W_k)/\G_{34}$ to a  model induced by
$f(W_k)\subset\tilde{N}/\D$. Notice that, by
Lemma~\ref{restrizionedimappe},
in the two processes we can choose the same $W_k$.\\
\underline{Step 6}: Combine Steps 4 and 5. More precisely, consider
the composition
$(f\circ\tau_{34})^{-1}\circ(f\circ\tau_{12})=\tau_{34}^{-1}\circ\tau_{12}$.
This is a diffeomorphism of the model $(\tau_{12}\circ
p_{12})^{-1}(W_k)/\G_{12}$ to the model $(\tau_{34}\circ
p_{34})^{-1}(W_k)/\G_{34}$. Such a diffeomorphism cannot exists, as
observed in Remark~\ref{impossibile}.\qed

\small
\appendix
\section{Appendix}
We now recall from the original article \cite{p1} the basic
definitions and results on quasifolds and related geometrical
objects. For some of them we give a reformulation that is suitable
for treating questions that arise in the study of diffeomorphisms
between quasifolds.

\begin{defn}[Quasifold model]\label{model}{\rm
Let $\ut$ be a $k$--dimensional connected smooth manifold and let
$\G$ be a discrete group acting by diffeomorphisms on $\ut$ so that
the set of points, $\ut_0$, where the action is not free, is closed
and has minimal codimension $\geq2$. This implies that the set
$\ut\setminus\ut_0$, where the action is free, is open, dense and
connected. Consider the space of orbits, $\ut/\G$, of the action of
the group $\G$ on $\ut$, endowed with the quotient topology, and the
canonical projection $p\;\colon\;\ut\rightarrow \ut/\G$. A {\em
quasifold model} of dimension $k$ is the triple $(\ut/\G,p,\ut)$,
shortly denoted $\ut/\G$.}
\end{defn}
\begin{defn}[Diffeomorphism of models]\label{diffeomodelli}{\rm
Given two models $(\ut/\G, p, \ut)$ and \linebreak $(\vt/\D,q,\vt)$,
a homeomorphism $f\,\colon\,\ut/\G\longrightarrow\vt/\D$ is a {\em
diffeomorphism} if there exists a diffeomorphism
$\ft\,\colon\,\ut\longrightarrow\vt$ such that $ q\circ \ft= f\circ
p$; we will then say that $\ft$ is a {\em lift} of $f$.}
\end{defn}
If the mapping $\ft$ is a lift of a diffeomorphism of models
$f\,\colon\,\ut/\G\longrightarrow \vt/\D$ so are the mappings
${\ft}^{\gamma}(-)=\ft(\gamma\cdot -)$, for all elements $\gamma$ in
$\G$, and $^{\delta}\ft(-)=\delta\cdot\ft(-)$, for all elements
$\delta$ in $\D$. We recall from \cite{p1} the following fundamental
Lemmas:
\begin{lemma}[Uniqueness of lifts]\label{orange}{\rm
Consider two models, $\ut/\G$ and $\vt/\D$, and let
$f\,\colon\,\ut/\G\longrightarrow\vt/\D$ be a diffeomorphism of
models. For any two lifts, $\ft$ and $\fbar$, of the diffeomorphism
$f$, there exists a unique element $\delta$ in $\D$ such that
$\fbar={}^\delta\ft$.}
\end{lemma}
\begin{lemma}[Equivariance of lifts]\label{green}{\rm Consider two models, $\ut/\G$ and $\vt/\D$, and
a diffeomorphism $f\,\colon\,\ut/\G\longrightarrow\vt/\D$. Then, for
a given lift, $\ft$, of the diffeomorphism $f$, there exists a group
isomorphism $F\,\colon\,\G\longrightarrow\D$ such that
${\ft}^{\gamma}={}^{F(\gamma)}\ft$, for all elements $\gamma$ in
$\G$.}
\end{lemma}
\begin{defn}[Smooth mapping between models]\label{cinftymodelli}{\rm
Given two models $(\ut/\G, p, \ut)$ and $(\vt/\D,q,\vt)$, a
continuous mapping $f\,\colon\,\ut/\G\longrightarrow\vt/\D$ is said
to be {\em smooth} if there exists a smooth mapping
$\ft\,\colon\,\ut\longrightarrow\vt$ and a homomorphism
$F\,\colon\,\G\rightarrow\D$ such that $\ft(\gamma
\tilde{u})=F(\gamma)\ft(\tilde{u})$; we will then say that $\ft$ is
an {\em equivariant lift} of $f$.}
\end{defn}
\begin{remark}[Induced model]\label{tuttistessopi}{\rm Let $(\ut/\G,p,\ut)$ be a
model. If $W$ is a connected open subset of the quotient $\ut/\G$,
then $p^{-1}(W)$ has countably many connected components; for any
two of them there is a $\gamma\in\G$ that takes the first one
diffeomorphically onto the second one. Let $\tilde{W}$ be a
connected component of $p^{-1}(W)$, let
$\G_{\tilde{W}}=\{\gamma\in\G\;|\;\gamma(\tilde{W})=\tilde{W}\}$ and
let $p_{\tilde{W}}=p_{|_{\tilde{W}}}$; then
$(\tilde{W}/\G_{\tilde{W}},p_{\tilde{W}},\G_{\tilde{W}})$ is a model
and $\tilde{W}/\G_{\tilde{W}}$ is homeomorphic to $W$. We will say
that $(\tilde{W}/\G_{\tilde{W}},p_{\tilde{W}},\G_{\tilde{W}})$ is a
{\em model induced} by $W\subset\ut/\G$. Notice that the models
induced by $W\subset\ut/\G$ are all diffeomorphic. }\end{remark} As
in \cite{p1}, the following definition is crucial for defining
quasifold structures:
\begin{defn}[Universal covering model]\label{simplyc}{\rm Consider a model of dimension
$k$, \linebreak $(\ut/\G,p,\ut)$. Let
$\pi\,\colon\,\ush\rightarrow\ut$ be the universal covering of $\ut$
and let $\hbox{Aut}(\ush|\ut)$ be the group of covering
automorphisms of $\ush\rightarrow\ut$. Then $\hbox{Aut}(\ush|\ut)$
acts on $\ush$ in a smooth, free and proper fashion with
$\ut=\ush/\hbox{Aut}(\ush|\ut)$. Consider the extension of the group
$\G$ by the group $\hbox{Aut}(\ush|\ut)$
$$1\longrightarrow \hbox{Aut}(\ush|\ut)\longrightarrow\Gsh
\longrightarrow\G\longrightarrow 1$$ defined as follows
$$\Gsh=\left\{\;\gsh\in\mbox{Diff}(\ush)\;|\;\exists\;
\gamma\in\Gamma\;\mbox{s. t.}\;\pi(\gsh(u^{\sharp}))=\gamma
\pi(u^{\sharp})\;\forall\; u^{\sharp}\in\ush\;\right\}.$$ It is easy
to verify that $\Gsh$ acts on $\ush$ according to the assumptions of
Definition~\ref{model} and that $\ut/\G=\ush/\Gsh$. If
$p^{\sharp}=p\circ\pi$ we will say that the triple
$(\ush/\Gsh,p^{\sharp},\ush)$ is the {\em universal covering model}
of $(\ut/\G,p,\ut)$. }\end{defn}
\begin{remark}[Uniqueness of the universal covering model]{\rm
If $(\ut/\G,p,\ut)$ is a model, then the uniqueness of the universal
covering implies that the universal covering model is unique up to
diffeomorphisms.}\end{remark}
\begin{remark}{\rm Let $(\ut/\G,p,\ut)$ be a
model and let $W$ be a connected open subset of the quotient
$\ut/\G$. For simplicity, in the sequel, instead of writing {\em
universal covering model of a model induced by} $W\subset\ut/\G$ we
will write {\em universal covering model induced by}
$W\subset\ut/\G$.}\end{remark}
\begin{prop}[Restriction property of $C^{\infty}$ mappings]\label{restrizionedimappe}{\rm Let
$f\,\colon\,\ut/\G\longrightarrow\vt/\D$ be a smooth mapping (a
diffeomorphism), let $W\subset\vt/\D$ be a connected open subset and
let $W_f$ be a connected component of $f^{-1}(W)$. Then the
restriction of $f\,\colon\, W_f\rightarrow W$ is a smooth mapping (a
diffeomorphism) from any model induced by $W_f\subset\ut/\G$ to a
corresponding model induced by $W\subset\vt/\D$; moreover, it is a
smooth mapping (a diffeomorphism) between the universal covering
models of these induced models.}\end{prop} \proof Let
$(\wt_f/\G_f,p_f,\wt_f)$ be a model induced by $W_f\subset\ut/\G$.
We can then choose a model $(\wt/\D_W,\D_W,\wt)$ such that
$\ft\,\colon\, \wt_f\rightarrow\wt$. Let
$F\,\colon\,\G\rightarrow\D$ the group homomorphism that makes $\ft$
equivariant. Notice that, if $\gamma\in\G_f$ and $w\in\wt_f$, then
$\ft(\gamma w)=F(\gamma)\ft(w)$, which in turn implies that
$F\,\colon\,\G_f\rightarrow\D_W$. Moreover, the restriction of $\ft$
to $\wt_f$ is a diffeomorphism if the mapping $f$ is a
diffeomorphism. This implies that the mapping $f$ is a smooth
mapping (diffeomorphism) between the induced models considered.
Consider the universal coverings
$\xymatrix{\wsh_f\ar[r]^{\pi_f}&\wt_f}$ and
$\xymatrix{\wsh\ar[r]^{\pi_W}&\wt}$. Since the fundamental group
$\pi_1(\wsh_f)$ is trivial, the mapping $\ft\circ\pi_f$ lifts to a
mapping $\xymatrix{\wsh_f\ar[r]^{\fsh}&\wsh}$ which is unique up to
the action of $\hbox{Aut}(\wsh|\wt)$. We obtain the following
diagram:
$$
\xymatrix@R=4pc@C=3pc{
W_f^{\sharp}\ar[r]^{\fsh}\ar[d]_{\pi_f}&W^{\sharp}\ar[d]^{\pi_W}\\
\wt_f\ar[d]_{p_f}\ar[r]^{\ft}&\wt\ar[d]^{p_W}\\
W_f\ar[r]^{f}&W}
$$
In order to prove the statement for smooth mappings we have to show
that there is a group homomorphism
$\xymatrix{\Gsh\ar[r]^{F^{\sharp}}&\dsh}$ such that $\fsh$ is
equivariant. Let $w$ be a point in $\wsh_f$ and let
$\gamma^{\sharp}\in\Gsh$; then
$\pi_W(\fsh(\gamma^{\sharp}w))=\ft\circ\pi_f(\gamma^{\sharp}w)=\ft(\gamma
\pi_f(w))=F(\gamma)\ft(\pi_f(w))$. Therefore there exists
$\delta^{\sharp}\in\dsh$ such that
$\fsh(\gamma^{\sharp}w)=\delta^{\sharp}\fsh(w)$ and
$\fsh\cdot\gamma^{\sharp}$ and $\delta^{\sharp}\cdot\fsh$ are both
lifts of $\ft\circ\pi_f$ that coincide at the point $w$. They are
thus equal on $\wsh$ by a basic result on coverings. This defines
the required homomorphism $\xymatrix{\Gsh\ar[r]^{F^{\sharp}}&\dsh}$.
If $\ft$ is a diffeomorphism, then $\fsh$ is a diffeomorphism and
$F^{\sharp}$ is an isomorphism. \qed
\begin{defn}[Differential form on a model]{\rm A {\em $k$--differential form}, $\omega$,
on a model \linebreak $(\vt/\G, p, \vt)$ is the assignment of a
$k$--differential form, $\widetilde{\omega}$, on $\ut$ that is
invariant by the action of the group $\G$.}
\end{defn}
\begin{defn}[Quasifold chart]\label{chart}{\rm
A dimension $k$ {\em quasifold chart} of a topological space $M$ is
the assignment of a quintuple $(U,\tau,\ut,p,\G)$ such that $U$ is a
connected open subset of $M$, $(\ut/\G, p, \ut)$ defines a
$k$--dimensional model, and the mapping $\tau$ is a homeomorphism of
the space $\ut/\G$ onto the set $U$. }\end{defn}
\begin{defn}[Compatible charts]\label{compatiblecharts}{\rm Two charts
$(\ua,\fia,\uta,\pa,\ga)$ and \linebreak $(\ub,\fib,\utb,\pb,\gb)$
such that $\ua\cap\ub\neq\emptyset$ are said to be {\em compatible}
if, for each $x\in\ua\cap\ub$, there exists an open connected subset
$U^x_{\alpha\beta}\subset\ua\cap\ub$ such that the homeomorphism
$$\gab=\fib^{-1}\circ\fia\;\colon\;\fia^{-1}(U^x_{\alpha\beta})\longrightarrow\fib^{-1}(U^x_{\alpha\beta})$$ is a
diffeomorphism from the universal covering model induced by
$\fia^{-1}(U^x_{\alpha\beta})\subset\uta/\ga$ to the universal
covering model induced by
$\fib^{-1}(U^x_{\alpha\beta})\subset\utb/\gb$. We will then say that
the mapping $\gab$ is a {\em change of charts}. }\end{defn}
\begin{defn}[Quasifold atlas]\label{atlas}
{\rm A dimension $k$ {\em quasifold atlas}, $\mathcal A$, on a
topological space $M$, is the assignment of a collection of
compatible charts
$${\mathcal A}= \{\,(\ua,\fia,\uta,\pa,\ga)\,|\,\alpha\in A\,\}$$ such that
the collection $\{\,\ua\,|\,\alpha\in A\,\}$ is an open cover of $M$.
}
\end{defn}
\begin{defn}[Complete atlas]\label{completeatlas}{\rm
The atlas $\cal A$ is {\em complete} if each chart compatible with
all of the charts in $\cal A$ belongs to $\cal A$.}
\end{defn}
\begin{prop}[Atlas completion]{\rm Let $M$ be a topological space endowed with a quasifold atlas $\mathcal A$.
Then there exists a unique complete atlas containing $\mathcal A$.}
\end{prop}
\proof This is a consequence of the two Lemmas that follow.\qed
\begin{lemma}[Restriction of a change of charts]\label{localizzazione}
{\rm Let $M$ be a quasifold, with quasifold structure given by the
atlas ${\cal A}=\{U_{\alpha},\alpha\in A\}$. Let $\alpha, \beta$ be
in $A$ such that $\ua\cap\ub\neq\emptyset$. Then, for each
$x\in\ua\cap\ub$ and for each connected open subset $W$ of
$U^x_{\alpha\beta}$, the homeomorphism
$$\gab=\fib^{-1}\circ\fia\;\colon\;\fia^{-1}(W)\longrightarrow\fib^{-1}(W)$$ is a
a diffeomorphism from the universal covering model induced by
$\fia^{-1}(W)\subset\fia^{-1}(U^x_{\alpha\beta})$ to the universal
covering model induced by
$\fib^{-1}(W)\subset\fib^{-1}(U^x_{\alpha\beta})$.}
\end{lemma}
\proof This is a consequence of
Proposition~\ref{restrizionedimappe}. We give the details in order
to establish some notation that will be useful later. By definition
of compatible charts we have the diffeomorphism of the universal
covering models:
\begin{equation}\label{cambiocarte}
\xymatrix@R=4pc@C=3pc{
(U^x_{\alpha\beta})^{\sharp}\ar[d]_{\piabx}\ar[rr]^{\gabx^{\sharp}}&&(U^x_{\beta\alpha})^{\sharp}\ar[d]^{\pibax}\\
\tilde{U}^x_{\alpha\beta}\ar[d]_{p_{\alpha,x}}&&\tilde{U}^x_{\beta\alpha}\ar[d]^{p_{\beta,x}}\\
\fia^{-1}(U^x_{\alpha\beta})\ar[r]^{\fia}\ar@/_2pc/[rr]_{\gab}&U^x_{\alpha\beta}\ar[r]^{\fib^{-1}}&\fib^{-1}(U^x_{\alpha\beta})}
\end{equation}
where
$(\tilde{U}_{\alpha\beta}^x/\G_{\alpha,x},p_{\alpha,x},\tilde{U}_{\alpha\beta}^x)$
is a model induced by $U^x_{\alpha\beta}\subset\ua$ and
$(\tilde{U}_{\beta\alpha}^x/\G_{\beta,x},p_{\beta,x},\tilde{U}_{\beta\alpha}^x)$
is a model induced by $U^x_{\alpha\beta}\subset\ub$. Consider a
model induced by $W\subset U^x_{\alpha\beta}\subset \ua$, with
$\tilde{W}_{\alpha}\subset \tilde{U}_{\alpha\beta}^x$:
$$(\tilde{W}_{\alpha}/\G_{\alpha,W},p_{\alpha,W},\tilde{W}_{\alpha}).$$
Now, refering to diagram~(\ref{cambiocarte}), apply
Proposition~\ref{restrizionedimappe} to the diffeomorphism $\gab$
and to the open subset
$\fib^{-1}(W)\subset\fib^{-1}(U^x_{\alpha\beta})$. More precisely,
we will have to consider the following items: a connected component,
$W^{\dagger}_{\alpha\beta}$, of $\piabx^{-1}(\tilde{W}_{\alpha})$;
the corresponding connected component, $W^{\dagger}_{\beta\alpha}$,
via $\gabx^{\sharp}$; the restriction, $\gabw^{\dagger}$, of
$\gabx^{\sharp}$ to $\tilde{W}^{\dagger}_{\alpha}$; the connected
component $\tilde{W}_{\beta}=\piabx(W^{\dagger}_{\beta})$ of
$\pb^{-1}(W)$ which yields a model,
$(\tilde{W}_{\beta}/\G_{\beta,W},p_{\beta,W},\tilde{W}_{\beta})$,
induced by $W\subset U^x_{\alpha\beta}\subset \ub$. We then get the
diagram
\begin{equation}\label{sollevamento}
\xymatrix@R=4pc@C=3pc{
W_{\alpha\beta}^{\sharp}\ar[d]\ar[rr]^{\gabw^{\sharp}}&&W_{\beta\alpha}^{\sharp}\ar[d]\\
W_{\alpha\beta}^{\dagger}\ar[d]\ar[rr]^{\gabw^{\dagger}}&&W_{\beta\alpha}^{\dagger}\ar[d]\\
\tilde{W}_{\alpha}\ar[d]_{p_{\alpha,W}}&&\tilde{W}_{\beta}\ar[d]^{p_{\beta,W}}\\
\fia^{-1}(W)\ar[r]^{\fia}\ar@/_2pc/[rr]_{\gab}&W\ar[r]^{\fib^{-1}}&\fib^{-1}(W)}
\end{equation}
where $W^{\sharp}_{\alpha\beta}\rightarrow W_{\alpha\beta}^{\dagger}$ and
$W^{\sharp}_{\beta\alpha}\rightarrow W_{\beta\alpha}^{\dagger}$ are universal coverings and the diffeomorphism
$\gabw^{\sharp}$ is the lift of $\gabw^{\dagger}$. This proves the lemma. \qed
\begin{lemma}[Mutual compatibility]\label{mutua}{\rm
Given a quasifold atlas $\mathcal A$ on a topological space $M$, we
have that, if $\ua$ and $\ub$ are charts that are compatible with
all of the charts of $\mathcal A$ and such that
$\ua\cap\ub\neq\emptyset$, then they are mutually compatible.}
\end{lemma}
\proof Consider $x\in\ua\cap\ub$; then there exists a chart $\uz$ in
$\mathcal A$ containing $x$. By assumption, this chart is compatible
with both $\ua$ and $\ub$. Let $U^x_{\alpha\zeta}$ and
$U^x_{\beta\zeta}$ be as in  Definition~\ref{atlas} and let us
define $U^x_{\alpha\beta}$ to be the connected component of
$U^x_{\alpha\zeta}\cap U^x_{\beta\zeta}$ containing $x$. To lighten
notation let us denote $U^x_{\alpha\beta}$ by $W$. By
Lemma~\ref{localizzazione}, the homeomorphisms
$\gaz\,\colon\,\fia^{-1}(W)\longrightarrow\fiz^{-1}(W)$ and
$\gzb\,\colon\,\fiz^{-1}(W)\longrightarrow\fib^{-1}(W)$ lift to
diffeomorphisms
$$\gazw^{\sharp}\,\colon\, W_{\alpha\zeta}^{\sharp}\longrightarrow W_{\zeta\alpha}^{\sharp}$$
and
$$\gzbw^{\sharp}\,\colon\, W_{\zeta\beta}^{\sharp}\longrightarrow W_{\beta\zeta}^{\sharp}.$$
Here $W^{\sharp}_{\zeta\alpha}\rightarrow
W^{\dagger}_{\zeta\alpha}\rightarrow \tilde{W}_{\zeta\alpha}$ and
$W^{\sharp}_{\zeta\beta}\rightarrow
W^{\dagger}_{\zeta\beta}\rightarrow \tilde{W}_{\zeta\beta}$ are
universal coverings of $\tilde{W}_{\zeta\alpha}$ and
$\tilde{W}_{\zeta\beta}$ respectively. Then, by
Remark~\ref{tuttistessopi}, there exists $\gamma\in\G_{\zeta,x}$
such that $\gamma(\tilde{W}_{\zeta\alpha})=\tilde{W}_{\zeta\beta}$.
If we consider its lift to the universal coverings we obtain the diagram:
\begin{equation}\label{zoomedcomposition}
\xymatrix@R=4pc@C=3pc{
\ar[r]^{\gazw^{\sharp}}&W^{\sharp}_{\zeta\alpha}\ar[r]^{\gamma^{\sharp}}\ar[d]&W^{\sharp}_{\zeta\beta}
\ar[r]^{\gzbw^{\sharp}}\ar[d]&\\
&W^{\dagger}_{\zeta\alpha}\ar[d]&    W^{\dagger}_{\zeta\beta}\ar[d]&\\
&\tilde{W}_{\zeta\alpha}\ar[r]^{\gamma}&\tilde{W}_{\zeta\beta}&
}
\end{equation}
Finally, we obtain a lift of the homeomorphism $\gab\,\colon\,
\fia^{-1}(W)\longrightarrow \fib^{-1}(W)$ by the following
composition of diffeomorphisms:
\begin{equation}\label{compo}
\gzbw^{\sharp}\circ \gamma^{\sharp}\circ\gazw^{\sharp}.
\end{equation}\qed
\begin{defn}[Quasifold structure]{\rm A {\em quasifold structure} on  a topological space $M$
is given by the assignment of a complete atlas.}\end{defn}
\begin{remark}\label{gruppino}{\rm
To each point $m\in M$ there corresponds a group $\G_m$ defined, up
to isomorphisms, as follows: given a chart $(\ua,\fia,\uta/\ga)$
around $m$, $\G_m$ is the isotropy group of $\ga$ at any point
$\tilde{u}\in \uta$ which projects down to $m$. The isomorphism
class of the group $\G_m$ does not depend  on the choice of the
point $\tilde{u}$ and of the chart. If all the $\G_m$'s are finite,
then $M$ is an orbifold; if they are trivial, then $M$ is a
manifold.}\end{remark}
\begin{prop}[Global quotient]\label{quozienteglobale}{\rm Let $\tilde{M}$ be an $n$--dimensional
connected smooth manifold, and let $\G$ be a discrete group acting
by diffeomorphisms on $\tilde{M}$ in such a way that the set of
points, $\tilde{M}_0$, where the action is not free, is closed and
has minimal codimension $\geq2$. The quotient $\tilde{M}/\G$ is an
$n$--dimensional quasifold. We will say that $\tilde{M}/\G$ is a
{\em global quotient}. }\end{prop} \proof The quotient
$M=\tilde{M}/\G$ is a model and is a quasifold covered by one chart.
We remark that, if $W\subset M$ is an open subset, then a model
induced by $W\subset M$ is obviously a chart compatible with
$\tilde{M}/\G$. Of course intersecting open subsets yields pairs of
compatible charts.\qed
\begin{defn}[Differential form]{\rm Let $M$ be a quasifold with quasifold structure
given by the atlas ${\cal A}=\{U_{\alpha},\alpha\in A\}$. A {\em
$k$--differential form}, $\omega$,  on $M$ is the assignment of a
$k$--differential form $\tilde{\omega}_{\alpha}$ on
$\tilde{U}_{\alpha}$ that is invariant by the action of
$\G_{\alpha}$, for each $\alpha\in A$. Moreover, the
$\tilde{\omega}_{\alpha}$'s must behave well under changes of
charts: if $\ua\cap\ub\neq\emptyset$ and $x\in \ua\cap\ub$, then
\begin{equation}\label{formecompatibili}
\piabx^*\tilde{\omega}_{\alpha}=
(\gabx^{\sharp})^*\circ
\pibax^*(\tilde{\omega}_{\beta})
\end{equation}
on
$(U^x_{\alpha\beta})^{\sharp}$.
}
\end{defn}
\begin{remark}[Properties of differential forms]\label{formedefinitesulcompletamento}{\rm Let $M$ be a
quasifold with structure defined by an atlas $\mathcal A$. Given a
form, $\omega$, on ${\mathcal A}$ the local forms are defined on all
of the charts of the completion of $\mathcal A$}\end{remark}
\begin{defn}[Diffeomorphism]\label{diffeo}{\rm A {\em diffeomorphism}, $f$, from a quasifold $M$,
defined by the atlas ${\cal A}$, to a quasifold $N$, defined by the
atlas ${\cal B}$, is given by a homeomorphism
$f\,\colon\,M\longrightarrow N$ with the following properties: for
each $x\in M$, there exist a chart $(\ua,\fia,\uta/\ga)$ containing
$x$, a chart $(V_b,k_b,\D_b,\tilde{V}_b)$ containing $f(x)$, and an
open subset $V_{\alpha b,x}$ of  $V_b$ such that $f(x)\in V_{\alpha
b,x}\subset V_b$, $f^{-1}(V_{\alpha b,x})\subset\ua$ and the
homeomorphism
$$k_b^{-1}\circ f\circ\fia\;\colon\;\fia^{-1}(f^{-1}(V_{\alpha b,x}))\rightarrow k_{b}^{-1}
(V_{\alpha b,x})$$ is a diffeomorphism from  the universal covering
model induced by $\fia^{-1}(f^{-1}(V_{\alpha b,x}))\subset\uta/\ga$
to the universal covering model induced by $k_{b}^{-1}(V_{\alpha
b,x})\subset\tilde{V}_b/\D_b$. The same holds for $f^{-1}$.}
\end{defn}
\begin{remark}\label{tuttelecarte}{\rm Let $f\,\colon\,M \rightarrow N$ be a diffeomorphism and let $x$ be a point in $M$.
The conditions in Definition~\ref{diffeo} are satisfied for each
pair of charts $U_{\alpha'}$ and $V_{b'}$ compatible with ${\mathcal
A}$ and ${\mathcal B}$ respectively and such that $x\in
U_{\alpha'}$, $f(x)\in V_{b'}$. This is an easy consequence of
Proposition~\ref{restrizionedimappe}: given a pair of charts
$U_{\alpha'}$ and $V_{b'}$ as above, we choose $V_{\alpha'b',x}$ to
be the connected component of the open subset $f(U_{\alpha'})\cap
V_{b'}\cap V_{\alpha b,x}$ that contains $f(x)$.}
\end{remark}
\begin{prop}[Diffeomorphism properties]\label{atlantevainatlante}{\rm If $f$ is a diffeomorphism from a quasifold $M$
defined by the atlas ${\cal A}$, to a quasifold $N$ defined by the
atlas ${\mathcal B}$, then, for each chart $(\ua,\fia,\uta/\ga)$
compatible with ${\mathcal A}$, the triple
$(f(\ua),f\circ\fia,\uta/\ga)$ is a chart in $N$ compatible with the
atlas $\mathcal B$. }
\end{prop}
\proof If $(\ua,\fia,\uta/\ga)\in {\mathcal A}$ is a chart, we need
to prove that the triple $(f(\ua),f\circ\fia,\uta/\ga)$ is a chart
compatible with the atlas $\mathcal B$. Let $(\vb,\kb,\vtb/\db)$ be
a chart in $\mathcal B$ such that $f(\ua)\cap\vb\neq\emptyset$, let
$y\in f(\ua)\cap\vb$ and consider and element $x\in\ua$ such that $f(x)=y$. Then, by
definition of diffeomorphism, there exists an open subset
$V_{\alpha\beta,x}$ in $\vb$ such that
$f^{-1}(V_{\alpha\beta,x})\subset\ua$ and the homeomorphism
$$k_{\beta}^{-1}\circ f\circ\fia\;\colon\;\fia^{-1}(f^{-1}(V_{\alpha\beta,x}))
\rightarrow k_{\beta}(V_{\alpha\beta,x})$$ is a diffeomorphism  of
the universal covering models of the respective induced models. This is exactly what is required for the two charts to be
compatible.\qed
\begin{defn}[Presmooth mapping]\label{presmoothmap}{\rm A {\em presmooth mapping}, $f$, from a quasifold $M$,
defined by the atlas ${\cal A}$, to a quasifold $N$, defined by the
atlas ${\cal B}$, is given by a continuous mapping
$f\,\colon\,M\longrightarrow N$ with the following properties: for
each $x\in M$, there exist a chart $(\ua,\fia,\uta/\ga)$ containing
$x$, a chart $(V_b,k_b,\D_b,\tilde{V}_b)$ containing $f(x)$ and an
open subset $V_{\alpha b,x}$ of  $V_b$ such that $f(x)\in V_{\alpha
b,x}\subset V_b$, $f^{-1}(V_{\alpha b,x})\subset\ua$ and the
continuous mapping
$$k_b^{-1}\circ f\circ\fia\;\colon\;\fia^{-1}(f^{-1}(V_{\alpha
b,x}))\rightarrow k_{b}^{-1} (V_{\alpha b,x})$$ is a
 smooth mapping from  the universal covering model induced by
$\fia^{-1}(f^{-1}(V_{\alpha b,x}))\subset\uta/\ga$ to the universal
covering model induced by $k_{b}^{-1}(V_{\alpha
b,x})\subset\tilde{V}_b/\D_b$.}
\end{defn}
\begin{remark}[From presmooth mapping to smooth mapping]\label{ultimo}{\rm We have already noticed some
differences between the definitions of smooth mappings and
diffeomorphisms in the local case. Let us see what we have to add to
the definition of presmooth mapping in order to obtain a smooth
mapping. Let $f$ be a presmooth mapping from a quasifold $M$,
defined by the atlas ${\cal A}$, to a quasifold $N$, defined by the
atlas ${\cal B}$. Take $x,x'\in M$, and, according to
Definition~\ref{presmoothmap}, take suitable charts
$(\ua,\fia,\uta/\ga)$ and $(\uapr,\fiapr,\utapr/\gapr)$ in $\mathcal
A$, suitable charts $(V_b,k_b,\D_b,\tilde{V}_b)$ and
$(V_b',k_b',\D_b',\tilde{V}_b')$ in $\mathcal B$, and $V_{\alpha
b,x}\subset V_b$, $V_{\alpha' b',x'}\subset V_b'$. Let $h$ be the
mapping $k_b^{-1}\circ f\circ\fia$ and $h'$ be the mapping
$k_{b'}^{-1}\circ f\circ\fiapr$. Suppose that $V_{\alpha b,x}\cap
V_{\alpha' b',x'}\neq\emptyset$, let $W$ be a connected component of
$V_{\alpha b,x}\cap V_{\alpha' b',x'}$ and let $W_f$ be a connected
component of $f^{-1}(W)$. Let us apply
Lemma~\ref{restrizionedimappe} to construct the diagram that follows
(here we use the notation we introduced in
diagram~(\ref{sollevamento})):
\begin{equation}\label{finale}
\xymatrix{
(W_f)^{\dagger}_{\alpha}\ar[d]\ar[rr]^{h^{\dagger}}&&W^{\dagger}_{b}\ar[d]\\
\fia^{-1}(W_f)\ar[dr]_{f\circ\fia}&&k^{-1}_b(W)\ar[ld]^{k^{-1}_b}\\
&W&\\
\fiapr^{-1}(W_f)\ar[ur]^{f\circ\fiapr}&&k^{-1}_{b'}(W)\ar[lu]_{k^{-1}_{b'}}\\
(W_f)^{\dagger}_{\alpha'}\ar[u]\ar[rr]^{(h')^{\dagger}}\ar@/^3pc/[uuuu]^{g_{\alpha\alpha',W_f}^{\dagger}}&&W^{\dagger}_{b'}\ar[u]\ar@/_3pc/[uuuu]_{g_{bb',W}^{\dagger}} \\
}\end{equation} If $f$ were
a diffeomorphism, by the uniqueness of the lift
(Lemma~\ref{orange}), we would have that there exists
$\delta^{\dagger}$ in the group acting on $W^{\dagger}_{b}$ such
that $h^{\dagger}={}^{\delta^{\dagger}}(g_{bb',W}^{\dagger}\circ
h'^{\dagger}\circ (g_{\alpha\alpha',W_f}^{\dagger})^{-1})$. In the
case of smooth mappings the uniqueness of the lift is no longer
guaranteed (an example in the orbifold case can be found in
\cite[Ex.~4.1.6b]{ch}). Therefore this gluing property has to be
required in the definition. It ensures, for example, that the
pullback of a form is well defined.}
\end{remark}
\begin{defn}[Smooth mapping]\label{smoothmap}{\rm
A {\em smooth mapping} $f$ from a quasifold $M$, defined by the
atlas ${\cal A}$, to a quasifold $N$, defined by the atlas ${\cal
B}$, is a presmooth mapping such that, for every $V_{\alpha b,x}$,
$V_{\alpha' b',x'}$ as in Remark~\ref{ultimo} and for each connected
open component, $W$, of $V_{\alpha b,x}\cap V_{\alpha' b',x'}$,
there exists $\delta^{\dagger}$ in the group acting on
$W^{\dagger}_{b}$ such that
$h^{\dagger}={}^{\delta^{\dagger}}(g_{bb',W}^{\dagger}\circ
h'^{\dagger}\circ (g_{\alpha\alpha',W_f}^{\dagger})^{-1})$. }
\end{defn}
\begin{remark}{\rm Let $M$ and $N$ be quasifolds, let $f\,\colon\,M\longrightarrow N$ be
a smooth mapping and let $\omega$ be a differential form on $N$.
Then the pullback $f^*(\omega)$ is a well defined form on $M$.}
\end{remark}

\noindent \sc Dipartimento di Matematica Applicata "G. Sansone",
Universit\`a di Firenze, Via S. Marta 3, 50139 Firenze,
ITALY, {\tt fiammetta.battaglia@unifi.it}\\
                        and\\
                        Dipartimento di Matematica "U. Dini", Universit\`a di Firenze,
                        Piazza Ghiberti 27, 50122 Firenze, ITALY, {\tt elisa.prato@unifi.it}
\end{document}

%% file: stargrid.pstex_t
\begin{picture}(0,0)%
\includegraphics{stargrid.eps}%
\end{picture}%
\setlength{\unitlength}{1697sp}%
\begingroup\makeatletter\ifx\SetFigFont\undefined%
\gdef\SetFigFont#1#2#3#4#5{%
  \reset@font\fontsize{#1}{#2pt}%
  \fontfamily{#3}\fontseries{#4}\fontshape{#5}%
  \selectfont}%
\fi\endgroup%
\begin{picture}(3705,3983)(3199,-6494)
\put(5601,-2674){\makebox(0,0)[lb]{\smash{{\SetFigFont{11}{13.2}{\rmdefault}{\bfdefault}{\updefault}{\color[rgb]{0,0,.82}$\rm{Y_1}$}%
}}}}
\put(6889,-4486){\makebox(0,0)[lb]{\smash{{\SetFigFont{11}{13.2}{\rmdefault}{\bfdefault}{\updefault}{\color[rgb]{0,0,.82}$\rm{Y_0}$}%
}}}}
\put(5589,-6424){\makebox(0,0)[lb]{\smash{{\SetFigFont{11}{13.2}{\rmdefault}{\bfdefault}{\updefault}{\color[rgb]{0,0,.82}$\rm{Y_4}$}%
}}}}
\put(3226,-5686){\makebox(0,0)[lb]{\smash{{\SetFigFont{11}{13.2}{\rmdefault}{\bfdefault}{\updefault}{\color[rgb]{0,0,.82}$\rm{Y_3}$}%
}}}}
\put(3214,-3349){\makebox(0,0)[lb]{\smash{{\SetFigFont{11}{13.2}{\rmdefault}{\bfdefault}{\updefault}{\color[rgb]{0,0,.82}$\rm{Y_2}$}%
}}}}
\end{picture}%

%% file: stargridduale.pstex_t
\begin{picture}(0,0)%
\includegraphics{stargridduale.eps}%
\end{picture}%
\setlength{\unitlength}{1697sp}%
\begingroup\makeatletter\ifx\SetFigFont\undefined%
\gdef\SetFigFont#1#2#3#4#5{%
  \reset@font\fontsize{#1}{#2pt}%
  \fontfamily{#3}\fontseries{#4}\fontshape{#5}%
  \selectfont}%
\fi\endgroup%
\begin{picture}(4043,3850)(2811,-6141)
\put(4816,-2446){\makebox(0,0)[lb]{\smash{{\SetFigFont{11}{13.2}{\familydefault}{\mddefault}{\updefault}{\color[rgb]{0,0,.82}$\rm{Y^*_0}$}%
}}}}
\put(6211,-6046){\makebox(0,0)[lb]{\smash{{\SetFigFont{11}{13.2}{\rmdefault}{\mddefault}{\updefault}{\color[rgb]{0,0,.82}$\rm{Y^*_3}$}%
}}}}
\put(3639,-6074){\makebox(0,0)[lb]{\smash{{\SetFigFont{11}{13.2}{\rmdefault}{\mddefault}{\updefault}{\color[rgb]{0,0,.82}$\rm{Y^*_2}$}%
}}}}
\put(2826,-3836){\makebox(0,0)[lb]{\smash{{\SetFigFont{11}{13.2}{\rmdefault}{\mddefault}{\updefault}{\color[rgb]{0,0,.82}$\rm{Y^*_1}$}%
}}}}
\put(6839,-3811){\makebox(0,0)[lb]{\smash{{\SetFigFont{11}{13.2}{\rmdefault}{\mddefault}{\updefault}{\color[rgb]{0,0,.69}$\rm{Y^*_4}$}%
}}}}
\end{picture}%